# A combined space discrete algorithm with a Taylor series by time for CFD

I.V. Kazachkov [*]

*Heat and Power Faculty, National Technical University of Ukraine "KPI", Kyiv, Polytekhnichna, 37, B.5, 02056, Ukraine*

**SUMMARY**

The first order by time partial differential equations are used as models in applications such as fluid flow, heat transfer, solid deformation, electromagnetic waves, and others. In this paper we propose the new numerical method to solve a class of initial-boundary value problems for the PDEs using one of the known space discrete numerical schemes and a Taylor series expansion by time. Normally a second order discretization by space is applied while a first order by time is satisfactory. Nevertheless, in a number of different problems, discretization both by temporal and by spatial variables is needed of highest orders, which complicates numerical solution, in some cases dramatically. Therefore it is difficult to apply the same numerical methods for the solution of some PDE arrays if their parameters are varying in a wide range so that in some of them different numerical schemes by time fit the best for precise numerical solution. The Taylor series based solution strategy for the non-stationary PDEs in CFD simulations has been proposed here that attempts to optimise the computation time and fidelity of the numerical solution.

## 1. INTRODUCTION

The second order PDEs have found extensive applications in the study of problems in fluid mechanics, flow in porous media, heat conduction, etc. [1-4]. A large number of numerical methods have been proposed for solving the second order PDEs, which are mainly the first order in time, in a CFD simulation. A key issue is the need to effectively use high performance numerical methods [5-12] and computers including parallel clusters [13] to complete analysis in time frames. In designing CFD software tools the author has attempted to build an essentially open single software framework, that enables arbitrarily complex non-stationary 3-D PDE array to be represented, which run efficiently on modern computers and allow simple increasing of accuracy in numerical simulation by time. The features of the above approach are that it employs Taylor series expansions to compute solution of the PDE by time without temporal discretization of the PDE. The idea of using the Taylor series expansions for numerical solution of non-stationary boundary problems has arisen from original use of a Taylor series described in [14,15] as an efficient procedure for parametric study in complex problems where a number of typical computations was replaced by Taylor series approximations.

The strategy for the numerical solution of non-stationary 3-D PDE using Taylor series by time has been proposed as follows:
• Numerical solution starts as usually with discretization of numerical domain and development of appropriate numerical grid
• Discretization of PDE by space is done by one of the known methods, which fits the best to the PDE given

---

[*] E-mail address: kazachkov@sotline.net.ua (I.V. Kazachkov).

- The temporal derivatives are computed up to the desired order by time for the numerical solution sought
- Using computed temporal derivatives the numerical solution sought is found from the Taylor series.

## 2. THE COMBINED NUMERICAL METHOD USING A TAYLOR SERIES BY TIME

Conventionally numerical solution of any initial-boundary problem for PDE or PDE array with one of the known numerical methods is going as follows:
1. Discretization of the numerical domain and development of the appropriate numerical grid.
2. Discretization of PDE by space and time with further transformation of the outgoing PDE to its approximation, for example algebraic finite-difference equations by space and time.
3. Numerical solution of the approximate (e.g. algebraic finite-difference) equations by space and time.
4. Testing of the numerical solution obtained and validation of it against the known data (other numerical solutions, analytical solutions for limit cases of the PDE stated, experimental results, etc.).

Highly important peculiarity of the above strategy is discretization of the PDE (step 2) performed according to the accuracy of the numerical solution by space and time required as far as this predetermines further steps and methods selected for the numerical solution. If any changes to the requirements of solution accuracy, then the step 2 changes, thus, the numerical algorithm changes totally. Our strategy replaces the steps 2, 3 of the above algorithm, which becomes the following one:
1. Spatial discretization of the numerical domain and development of the appropriate numerical grid.
2. Discretization of PDE by space with further transformation of the outgoing PDE to its approximation, for example, algebraic finite-difference equations by space.
3. Numerical solution of the approximate (e.g. algebraic finite-difference) equations by space.
4. Computing the temporal derivatives using the outgoing PDE or PDE array with further calculation of the numerical solution in time based on the Taylor series expansion by time.
5. Testing the numerical solution obtained and its validation against the known data (other numerical solutions, analytical solutions for limit cases of the PDE stated, experimental results, etc.).

Let us start with a few simple examples showing the idea of the proposed strategy. For this, first consider the following one-dimensional non-stationary equation (describing, for example 1-D flow) along with the corresponding initial and boundary conditions:

$$\frac{\partial U}{\partial t} + U \frac{\partial U}{\partial x} = \frac{\partial^2 U}{\partial x^2} \tag{1}$$

$$t = 0, \quad U = U_0(x), \quad x \in \Gamma, \ U = U_\Gamma(t). \tag{2}$$

We do not specify the boundary condition (2) because it is no matter for explanation of the proposed method.

Supposed $\frac{\partial^2 U}{\partial x \partial t} = \frac{\partial^2 U}{\partial t \partial x}$, the equation (1) is rewritten in a more convenient form

$$\frac{\partial U}{\partial t} = -U\frac{\partial U}{\partial x} + \frac{\partial^2 U}{\partial x^2} = f(x), \tag{3}$$

where from, differentiating the last equation by time, results

$$\frac{\partial^2 U}{\partial t^2} = \frac{\partial}{\partial t}\left(\frac{\partial^2 U}{\partial x^2} - U\frac{\partial U}{\partial x}\right) = \frac{\partial^2}{\partial x^2}\left(\frac{\partial U}{\partial t}\right) - \frac{\partial U}{\partial t}\frac{\partial U}{\partial x} - U\frac{\partial}{\partial x}\left(\frac{\partial U}{\partial t}\right) = \frac{\partial^2 f}{\partial x^2} - \frac{\partial}{\partial x}(fU). \tag{4}$$

Now a Taylor series with a second order accuracy by time (introduce $\Delta t$ as the time step in numerical solution), accounted (3), (4) yields the following approximation of the solution to the equation (1)

$$U = U_0 + f_0 \Delta t + \frac{1}{2!}\left(\frac{\partial^2 f}{\partial x^2} - \frac{\partial(fu)}{\partial x}\right)_0 (\Delta t)^2 + o\left((\Delta t)^2\right), \tag{5}$$

Where $U_0$ is the known initial data, and

$$f_0 = \left(\frac{\partial^2 U}{\partial t^2}\right)_0 - U_0\left(\frac{\partial U}{\partial x}\right) = \frac{\partial^2 U_0}{\partial x^2} - U_0\frac{\partial U_0}{\partial x} + o\left((\Delta t)^2\right), \tag{6}$$

is easily computed by the known function $U_0(x)$. In a similar way,

$$\left(\frac{\partial^2 f}{\partial x^2} - \frac{\partial(fU)}{\partial x}\right)_0 = \frac{\partial^4 U_0}{\partial x^4} - \frac{\partial^2}{\partial x^2}\left(U_0\frac{\partial U_0}{\partial x}\right) - \frac{\partial}{\partial x}\left(U_0\frac{\partial^2 U_0}{\partial x^2} - U_0^2\frac{\partial U_0}{\partial x}\right). \tag{7}$$

Then substitution of the equations (6), (7) into a Taylor series (5) results in

$$U = U_0 + \left(\frac{\partial^2 U}{\partial x^2} - U\frac{\partial U}{\partial x}\right)_0 \Delta t + \frac{1}{2}\left[\frac{\partial^4 U_0}{\partial x^4} - 2U_0\frac{\partial^3 U_0}{\partial x^3} + \left(U_0^2 - 4\frac{\partial U_0}{\partial x}\right)\frac{\partial^2 U_0}{\partial x^2}\right.$$
$$\left. +2U_0\left(\frac{\partial U_0}{\partial x}\right)^2\right](\Delta t)^2 + o\left((\Delta t)^2\right), \tag{8}$$

and so on. Evidently, one can continue this procedure to get any desired order of accuracy by time. The transformation of the procedure from any n–th layer by time to the (n+1)-th layer by time is similar to the stated in the equations (5), (8) above:

$$U_{n+1} = U_n + f_n \Delta t + \frac{1}{2}\left(\frac{\partial^2 f}{\partial x^2} - \frac{\partial(fU)}{\partial x}\right)_n (\Delta t)^2 + o\left((\Delta t)^2\right), \tag{9}$$

where n=0,1,2,…,N.

Thus, the right hand of the equations (5), (9) is always a function of the coordinates at the current moment of time. Thus, neither explicit, nor implicit approximations by time are applied; no difference equations by time are needed! As the equation (9) shows, the second order by time approximation in the equation (1) requests all spatial derivatives of the function sought up to the 4-th order. Adding the next term in the time series requires the corresponding term twice differentiated by space. If computing the highest order derivatives $\partial^i f / \partial t^i$ is analytically complicated, it is done numerically. Consequently, instead of a solution of a difference (or any discrete) equation, computation of the spatial derivatives with the accuracy stated is proposed. Then the numerical solution sought is computed from the Taylor series by time.

### 3. EXAMPLES OF A TAYLOR SERIES IN NUMERICAL SOLUTIONS

**Example 1**. In case of the simple wave equation with the following initial and boundary conditions

$$\frac{\partial u}{\partial t} = \frac{\partial u}{\partial x}; \qquad t = 0, \quad u = u_0(x) = x; \qquad x=0, \quad u = u(t) = t \tag{10}$$

an analytical solution of the boundary problem (10) is known as a wave spreading with the velocity 1 countercurrent to the axis $x$, $u = f(x+t)$. According to our strategy, the numerical solution of (10) is $u = u_0 + (\partial u / \partial t)_0 \Delta t + O((\Delta t)^2)$, where from with account of the above-mentioned $u = x + \Delta t$. The solution obtained does not change with an increase of accuracy because the first order solution coincides here with the exact analytical solution.

**Example 2.** The one-dimensional non-linear equation

$$\frac{\partial U}{\partial t} + U \frac{\partial U}{\partial x} = v \frac{\partial^2 U}{\partial x^2}, \tag{11}$$

with the following initial and boundary conditions

$$t = 0, \quad U = 1; \qquad x=0, \quad \frac{\partial U}{\partial x} = 0, \quad \frac{\partial^2 U}{\partial x^2} = 0, \tag{12}$$

is solved according to the proposed strategy as follows

$$\frac{\partial U}{\partial t} = v \frac{\partial^2 U}{\partial x^2} - U \frac{\partial U}{\partial x},$$

$$\frac{\partial^2 U}{\partial t^2} = \frac{\partial}{\partial t}\left(\frac{\partial U}{\partial t}\right) = v \frac{\partial^2}{\partial x^2}\left(\frac{\partial U}{\partial t}\right) - \frac{1}{2}\frac{\partial}{\partial t}\frac{\partial U^2}{\partial x} = v \frac{\partial^2}{\partial x^2}\left(\frac{\partial U}{\partial t}\right) - \frac{\partial}{\partial x}\left(U \frac{\partial U}{\partial t}\right) =$$

$$= v \frac{\partial^2}{\partial x^2}\left(v \frac{\partial^2 U}{\partial x^2} - U \frac{\partial U}{\partial x}\right) - \frac{\partial}{\partial x}\left[U\left(v \frac{\partial^2 U}{\partial x^2} - U \frac{\partial U}{\partial x}\right)\right],$$

where from yields solution to the initial-boundary problem (11), (12) in a second order approach by time

$$U = U_0 + \left(v \frac{\partial^2 U}{\partial x^2} - U \frac{\partial U}{\partial x}\right)\Delta t + \frac{1}{2}\left\{v \frac{\partial^2}{\partial x^2}\left(v \frac{\partial^2 U}{\partial x^2} - U \frac{\partial U}{\partial x}\right) - \frac{\partial}{\partial x}\left[U\left(v \frac{\partial^2 U}{\partial x^2} - U \frac{\partial U}{\partial x}\right)\right]\right\}(\Delta t)^2$$

The first order approach by time gives here the same result as any higher order approaches: $U = 1$.

## 4. ALGORITHM FOR NUMERICAL SOLUTION OF THE NAVIER-STOKES EQUATIONS WITH THE PROPOSED STRATEGY

Consider numerical algorithm for the solution to the equations of 3-D non-stationary motion of heat-conducting incompressible viscous fluids:

$$div\,\vec{v}=0,\quad \frac{\partial \vec{v}}{\partial t}+\vec{v}\nabla\vec{v}=-\frac{1}{\rho}\nabla p+\nu\Delta\vec{v}+\vec{f},\quad \frac{\partial T}{\partial t}+\vec{v}\nabla T=\frac{1}{\rho c}div(\lambda\nabla T)+\Phi, \qquad (13)$$

where $\vec{v}=\{v_x,v_y,v_z\}$, $\vec{f}=\{f_x,f_y,f_z\}$ are the velocity and the external force vectors, respectively. The Cartesian coordinates $x,y,z$ are implied here, then $\rho,\mu,\lambda$ are the density, dynamic viscosity coefficient, and heat conductivity coefficient, correspondingly, $\nu$ the coefficient of kinematical viscosity and $c$ the specific heat capacity, $\nabla,\Delta$ denote the gradient and Laplace operators. Finally, $\Phi$ is the dissipative function,

$$\Phi=\frac{\mu}{\rho c}\left(2|\nabla\vec{v}|^2+\left|\frac{\tau}{\mu}\right|^2\right),\quad \left|\frac{\tau}{\mu}\right|^2=\left(\frac{\partial v_x}{\partial y}+\frac{\partial v_y}{\partial x}\right)^2+\left(\frac{\partial v_x}{\partial z}+\frac{\partial v_z}{\partial x}\right)^2+\left(\frac{\partial v_y}{\partial z}+\frac{\partial v_z}{\partial y}\right)^2.$$

The partial differential equation array (13) thus obtained has to be supplemented with the corresponding initial conditions stated in the numerical domain $\Omega$:

$$t=0,\quad \vec{v}=\vec{v}_0(x,y,z),\quad p=p_0(x,y,z),\quad T=T_0(x,y,z),\quad (x,y,z)\in\Omega, \qquad (14)$$

as well as with the corresponding boundary conditions (Dirichlet, Neumann, mixed, etc.) at the boundary $(x,y,z)\in\Gamma$. They are not specified here because it has no matter for the proposed strategy of the numerical solution, which is applicable by any boundary conditions. To apply the above described numerical strategy to the Navier-Stokes equation array (13) with the initial conditions (14), rewrite these equations in the form:

$$div\,\vec{v}=0,\quad \frac{\partial\vec{v}}{\partial t}=\vec{F}(\vec{v},\nabla\vec{v},\vec{f})-\frac{1}{\rho}\nabla p,\quad \frac{\partial T}{\partial t}=F_T(\vec{v},T,\nabla T,\nabla\vec{v}), \qquad (15)$$

where is

$$\vec{F}(\vec{v},\nabla\vec{v},\vec{f})=\nu\Delta\vec{v}-\vec{v}\nabla\vec{v}+\vec{f},\quad F_T(\vec{v},T,\nabla T,\nabla\vec{v})=\Phi+\frac{1}{\rho c}div(\lambda\nabla T)-\vec{v}\nabla T, \qquad (16)$$

$\vec{F}=(F_x,F_y,F_z)$ is the vector of the right hand of the momentum equation except the pressure gradient.

Thus, all right hands of the equations (16) are known at the initial moment of time from the initial data (14), and then all temporal derivatives of the velocity vector and of the temperature are computed from (15). Consequently, the velocity and temperature fields are calculated at the next time step from the Taylor series:

$$\vec{v}=\vec{v}_0(x)+\left(\frac{\partial\vec{v}}{\partial t}\right)_0\Delta t+\frac{1}{2!}\left(\frac{\partial^2\vec{v}}{\partial t^2}\right)_0(\Delta t)^2+\frac{1}{3!}\left(\frac{\partial^3\vec{v}}{\partial t^3}\right)_0(\Delta t)^3+o((\Delta t)^3), \qquad (17)$$

$$T = T_0(x) + \left(\frac{\partial T}{\partial t}\right)_0 \Delta t + \frac{1}{2!}\left(\frac{\partial^2 T}{\partial t^2}\right)_0 (\Delta t)^2 + \frac{1}{3!}\left(\frac{\partial^3 T}{\partial t^3}\right)_0 (\Delta t)^3 + o((\Delta t)^3).$$

Here $\Delta t$ is the temporal step chosen for the generated numerical grid. The approximate numerical solution (17) is computed with a required order of accuracy by time (here it is up to the third order terms for example). It is very important that at the first temporal step computed by equations (17) the pressure distribution is unknown and continuity equation has not been used yet. Surprisingly, the velocity and the temperature fields do not depend on the pressure distribution at the first time step as shown in detail below.

This numerical scheme in a first order by $\Delta t$ (when only the terms up to $\Delta t$ are kept in (17)) completely coincide with the simplest first order explicit numerical scheme. But these two methods differ a lot afterwards. For example, well-known numerical schemes of a second order by time are very time-consuming and cumbersome while in the strategy proposed here the numerical solution procedure in a second order by time (as well as in any higher order by time) is nearly the same as in the first order by time. Moreover, any highest order numerical solution is got similarly and, what is very important, does not request more computer resources than the first order solution. All needed for this is just an easy computation by equations (15), (16) with further substitution of the results into the Taylor series (17).

The numerical solution of a second order accuracy by time requests, in contrast with the first order approximation, calculation of the derivatives of pressure, because the equation array (13) transforms to

$$\frac{\partial}{\partial t}(div\,\vec{v}) = 0, \quad \frac{\partial}{\partial t}\left(\frac{\partial \vec{v}}{\partial t}\right) = \frac{\partial^2 \vec{v}}{\partial t^2} = \frac{\partial \vec{F}}{\partial t} - \frac{1}{\rho}\frac{\partial}{\partial t}\nabla p, \quad \frac{\partial^2 T}{\partial t^2} = \frac{\partial}{\partial t}\left(\frac{\partial T}{\partial t}\right) = \frac{\partial F_T}{\partial t}, \quad (18)$$

where

$$\frac{\partial \vec{F}}{\partial t} = \nu\Delta\left(\frac{\partial \vec{v}}{\partial t}\right) + \frac{\partial \vec{f}}{\partial t} - \left[\left(\frac{\partial \vec{v}}{\partial t}\right)\nabla\vec{v} + \vec{v}\nabla\left(\frac{\partial \vec{v}}{\partial t}\right)\right],$$

$$\frac{\partial F_T}{\partial t} = \Phi_1 + \frac{1}{\rho c}div\left[\lambda\nabla\left(\frac{\partial T}{\partial t}\right)\right] - \vec{v}\nabla\left(\frac{\partial T}{\partial t}\right) - \nabla T\left(\frac{\partial \vec{v}}{\partial t}\right).$$

Here

$$\Phi_1 = \frac{\partial \Phi}{\partial t} = \frac{\mu}{\rho c}\left(2\left|\nabla\vec{F}\right|^2 + \left|\frac{\tau_1}{\mu}\right|^2\right),$$

$$\left|\frac{\tau_1}{\mu}\right|^2 = \left(\frac{\partial F_x}{\partial y} + \frac{\partial F_y}{\partial x}\right)^2 + \left(\frac{\partial F_x}{\partial z} + \frac{\partial F_z}{\partial x}\right)^2 + \left(\frac{\partial F_y}{\partial z} + \frac{\partial F_z}{\partial y}\right)^2.$$

Then further transformation results in

$$\vec{F}_1 = \frac{\partial \vec{F}}{\partial t} = \nu\Delta\left(\vec{F} - \frac{1}{\rho}\nabla p\right) + \frac{\partial \vec{f}}{\partial t} - \left[\left(\vec{F} - \frac{1}{\rho}\nabla p\right)\nabla\vec{v} + \vec{v}\left(\nabla\vec{F} - \frac{1}{\rho}\Delta p\right)\right],$$

$$F_T^1 = \frac{\partial F_T}{\partial t} = \Phi_1 + \frac{1}{\rho c}div(\lambda\nabla F_T) - \left[\vec{v}\nabla F_T + \nabla T\left(\vec{F} - \frac{1}{\rho}\nabla p\right)\right], \quad (19)$$

And now transformation of the first equation in the equation array (18) to the form $div(\partial \vec{v}/\partial t)=0$, with further substitution of the components $(\partial \vec{v}/\partial t)$ from the equations (15), yields

$$div(\partial \vec{v}/\partial t) = div\vec{F} - (1/\rho)div(\nabla p) = 0,$$

where from finally goes to the following Poisson equation

$$\Delta p = \rho\, div\vec{F}\,. \tag{20}$$

The equation (20) thus obtained allows computing the pressure distribution in the numerical domain by the known values of the vector $\vec{F}$, which have been computed from the equation (16) based on the first order approximation for the velocity field (described above). This equation (20) is solved comparably easily and further numerical algorithm does not request solution of any equation because all needed here is computing the spatial derivatives for the functions in the numerical domain. The approximation of the second order accuracy by time to the numerical solution is got afterwards simply from the Taylor series (17).

Finally, the closed system of the equations (18)-(20) has been got here for the computation of the second order by time numerical approximation to the outgoing Navier-Stokes equation. Importantly, any difference equations are absent here, except the one well-known Poisson equation (20) for the pressure. Only the spatial derivatives are to be computed with the required accuracy, which is much easier than solution of the difference equations and does not complicate the algorithm in a second order approach.

Now differentiating the partial differential equations (18) by time, accounting (19), the third order accuracy by time is got as follows:

$$\frac{\partial^3 \vec{v}}{\partial t^3} = \frac{\partial^2 \vec{F}}{\partial t^2} - \frac{1}{\rho}\frac{\partial^2}{\partial t^2}\nabla p, \qquad \frac{\partial^3 T}{\partial t^3} = \frac{\partial^2 F_T}{\partial t^2}, \tag{21}$$

where $\dfrac{\partial^2}{\partial t^2}\nabla p$ is computed based on the first and on the second approximations by time got for the pressure and its derivative by time, including the initial data as described above. Obviously, the second order derivative by time for the pressure can be computed only with the first order accuracy because the first order derivative has been computed with the first order accuracy and afterwards the second order derivatives were computed also with the first order accuracy as derivative from the first order derivative. This important question is subject for a separate detail investigation. Here are

$$\frac{\partial^2 \vec{F}}{\partial t^2} = \frac{\partial}{\partial t}\left(\frac{\partial \vec{F}}{\partial t}\right) = \nu\Delta\left[\left(\frac{\partial \vec{F}}{\partial t}\right) - \frac{1}{\rho}\frac{\partial}{\partial t}\nabla p\right] + \frac{\partial^2 \vec{f}}{\partial t^2} - \left\{\left[\left(\frac{\partial \vec{F}}{\partial t}\right) - \frac{1}{\rho}\frac{\partial}{\partial t}\nabla p\right]\nabla \vec{v} + \right.$$

$$\left. + \left(\vec{F} - \frac{1}{\rho}\nabla p\right)\nabla\left(\vec{F} - \frac{1}{\rho}\nabla p\right) + \left(\vec{F} - \frac{1}{\rho}\nabla p\right)\nabla\left(\vec{F} - \frac{1}{\rho}\nabla p\right) + \vec{v}\left[\nabla\left(\frac{\partial \vec{F}}{\partial t}\right) - \frac{1}{\rho}\Delta\left(\frac{\partial p}{\partial t}\right)\right]\right\} =$$

$$= \nu\Delta\left[\left(\frac{\partial \vec{F}}{\partial t}\right) - \frac{1}{\rho}\frac{\partial \nabla p}{\partial t}\right] + \frac{\partial^2 \vec{f}}{\partial t^2} - \left\{\left[\left(\frac{\partial \vec{F}}{\partial t}\right) - \frac{1}{\rho}\frac{\partial \nabla p}{\partial t}\right]\nabla \vec{v} + \right.$$

$$\left. + 2\left(\vec{F} - \frac{\nabla p}{\rho}\right)\left(\nabla \vec{F} - \frac{\Delta p}{\rho}\right) + \vec{v}\left[\nabla\left(\frac{\partial \vec{F}}{\partial t}\right) - \frac{1}{\rho}\Delta\left(\frac{\partial p}{\partial t}\right)\right]\right\},$$

$$\frac{\partial^2 F_T}{\partial t^2} = \frac{\partial}{\partial t}\left(\frac{\partial F_T}{\partial t}\right) = \Phi_2 + \frac{1}{\rho c} div\left(\lambda \nabla F_T^1\right) - \left(\vec{F}_1 \nabla T + 2\vec{F}\nabla F_T + \vec{v}\nabla F_T^1\right),$$

where

$$\Phi_2 = \frac{\partial \Phi_1}{\partial t} = \frac{\mu}{\rho c}\left(2\left|\nabla \vec{F}_1\right|^2 + \left|\frac{\tau_2}{\mu}\right|^2\right),$$

$$\left|\frac{\tau_2}{\mu}\right|^2 = \left(\frac{\partial F_x^1}{\partial y}+\frac{\partial F_y^1}{\partial x}\right)^2 + \left(\frac{\partial F_x^1}{\partial z}+\frac{\partial F_z^1}{\partial x}\right)^2 + \left(\frac{\partial F_y^1}{\partial z}+\frac{\partial F_z^1}{\partial y}\right)^2.$$

$\vec{F}_1 = \left(F_x^1, F_y^1, F_z^1\right)$, $F_T^1$ are the new right hands of the equations according to (19).

Finally the equations for the solution of a third order accuracy by time (21) result in

$$\frac{\partial^3 \vec{v}}{\partial t^3} = \vec{F}_2 - \frac{1}{\rho}\frac{\partial^2}{\partial t^2}\nabla p, \qquad \frac{\partial^3 T}{\partial t^3} = F_T^2,$$

$$F_T^2 = \Phi_2 + \frac{1}{\rho c}div\left(\lambda \nabla F_T^1\right) - \left(\vec{F}_1\nabla T + 2\vec{F}\nabla F_T + \vec{v}\nabla F_T^1\right),$$

$$\vec{F}_2 = \frac{\partial^2 \vec{f}}{\partial t^2} - \left\{\left[\vec{F}_1 - \frac{1}{\rho}\frac{\partial \nabla p}{\partial t}\right]\nabla \vec{v} + 2\left(\vec{F} - \frac{\nabla p}{\rho}\right)\left(\nabla \vec{F} - \frac{\Delta p}{\rho}\right) + \vec{v}\left[\nabla \vec{F}_1 - \frac{1}{\rho}\Delta\left(\frac{\partial p}{\partial t}\right)\right]\right\}$$

$$+\nu \Delta\left(\vec{F}_1 - \frac{1}{\rho}\frac{\partial \nabla p}{\partial t}\right). \tag{22}$$

Here $\partial^2 \vec{f}/\partial t^2$ is known as a derivative from an external force stated. The most often it is constant (e.g. gravitation) or known (stated electromagnetic force in a conductive media, acceleration due to vibration, etc).

Consequently, the equations (22) give the numerical solution for the outgoing equations (15), which is a third order accuracy by time. The pressure distribution in a third order accuracy is got from solution of the Poisson equation (20) after substitution of the computed velocity field. Obviously, there is no problem to implement the proposed method to the case of variable physical properties of continua and to some other more general cases. It must be noted that an increase of accuracy of the approximations by time requests computing the temporal derivatives of the pressure, which needs to keep in computer memory a few temporal layers. Therefore this question needs a separate deep study. Nevertheless, this preliminary analysis shown high efficiency and simplicity of the strategy proposed for the solution of the non-stationary non-isothermal equations of the Navier-Stokes type, as well as many other first orders by time PDE arrays.

The same strategy is applied consequently at each and every time step and the transformation from any n-th layer by time to the (n+ 1)-th layer by time ($t_{n+1} = t_n + \Delta t$) is similar to the stated in the equations (17):

$$\vec{v}_{n+1} = \vec{v}_n(x) + \left(\frac{\partial \vec{v}}{\partial t}\right)_n \Delta t + \frac{1}{2!}\left(\frac{\partial^2 \vec{v}}{\partial t^2}\right)_n (\Delta t)^2 + \frac{1}{3!}\left(\frac{\partial^3 \vec{v}}{\partial t^3}\right)_n (\Delta t)^3 + o((\Delta t)^3),$$

$$T_{n+1} = T_n(x) + \left(\frac{\partial T}{\partial t}\right)_n \Delta t + \frac{1}{2!}\left(\frac{\partial^2 T}{\partial t^2}\right)_n (\Delta t)^2 + \frac{1}{3!}\left(\frac{\partial^3 T}{\partial t^3}\right)_n (\Delta t)^3 + o((\Delta t)^3). \tag{23}$$

where n=0,1,2,…,N is. Certainly, the time steps $\Delta t$ and the number of iterations at each time step may vary from one time step to another time step; therefore there is also the separate subject for additional investigation concerning the optimization of the numerical algorithm.

## 5. COMPARISON OF THE PROPOSED STRATEGY WITH THE METHOD OF FRACTIONAL DIFFERENTIALS

Let us consider an example from [16] on computation of a heat flux by the known temperature distribution using an analytical method of the fractional differentials. The task on heating of the semi-infinite domain is modeled by the following equation with the corresponding initial and boundary conditions:

$$\left(\frac{\partial}{\partial t} - \frac{\partial^2}{\partial x^2}\right)T = 0, \quad 0 < x < \infty, \quad 0 < t < \infty, \quad T_{x=0} = T_s(t), \quad T_{x=\infty} = 0, \quad T_{t=0} = 0. \tag{24}$$

The heat flux is $q_s = (\partial T / \partial x)_{x=0}$. The differential operator in (24) may be represented in the form [16]:

$$\left(\frac{\partial^{1/2}}{\partial t^{1/2}} - \frac{\partial}{\partial x}\right)\left(\frac{\partial^{1/2}}{\partial t^{1/2}} + \frac{\partial}{\partial x}\right)T = 0, \tag{25}$$

supposed that $\dfrac{\partial^{1/2}T}{\partial t^{1/2}}\dfrac{\partial T}{\partial x} = \dfrac{\partial}{\partial x}\dfrac{\partial^{1/2}T}{\partial t^{1/2}}$. Now consider the equation presented by the right multiplayer of (25):

$$\left(\frac{\partial^{1/2}}{\partial t^{1/2}} + \frac{\partial}{\partial x}\right)T = 0. \tag{26}$$

The solutions to the equation (26) are also solutions to the equations (25) because the operator applied to a zero results in zero. Thus, solutions to (26) are solutions to (24) as well. The equation (4) written for $x=0$ gives immediately solution of the task stated, namely temperature gradient at the boundary of the domain:

$$-q_s = -\left(\frac{\partial T}{\partial x}\right)_0 = \frac{d^{1/2}T_s(t)}{dt^{1/2}} = \frac{1}{\sqrt{\pi}}\frac{d}{dt}\int_0^t \frac{T_s(t)}{\sqrt{t-\tau}}d\tau. \tag{27}$$

Note that the temperature gradient (27) has been found without solution of the task for temperature distribution (24). This is why the method of the fractional differentials is also called the non-field method. It allows computing analytically a heat flux at the boundary of a domain directly through such comparably simple transformation of the outgoing differential equation.

Let us analyze the simple linear equation to compare the method of fractional differentials with the method proposed here. For this, except solution (27), consider also a general solution to the boundary problem (24) [17]:

$$T(x,t) = -\int_0^t T_s(\tau)\frac{\partial T_\tau}{\partial \tau}d\tau, \tag{28}$$

$$-\frac{\partial T_\tau}{\partial \tau} = \frac{\partial}{\partial \tau}\Theta\left(\frac{x}{2\sqrt{t-\tau}}\right) = \frac{\partial}{\partial \tau}\frac{2}{\sqrt{\pi}}\int_0^{\frac{x}{2\sqrt{t-\tau}}} e^{-\alpha^2}d\alpha = \frac{x}{2\sqrt{\pi}(t-\tau)^{3/2}}e^{-\frac{x^2}{4(t-\tau)}}, \quad (29)$$

the temperature $T_s(t)$ is kept by $x=0$ for all $\tau$ from 0 till $t$. From (28), (29) follows

$$T(x,t) = \frac{x}{2\sqrt{\pi}}\int_0^t \frac{T_s(\tau)}{(t-\tau)^{3/2}}e^{-\frac{x^2}{4(t-\tau)}}d\tau. \quad (30)$$

Introducing in (30) the new variable $\xi = \frac{x}{2\sqrt{t-\tau}}$ results in

$$T(x,t) = \frac{2}{\sqrt{\pi}}\int_{\frac{x}{2\sqrt{t}}}^\infty T_s\left(t+\frac{x^2}{4\xi^2}\right)e^{-\xi^2}d\xi \quad (31)$$

which gives solution to the task (24). By $x=0$, the solution (31) satisfies to the boundary condition (24).

Now $q_s$ can be got from (31) or directly from (30) with differentiation of the integral by parameter [18]:

$$\frac{\partial T}{\partial x} = \frac{1}{2\sqrt{\pi}}\int_0^t \frac{T_s(\tau)}{(t-\tau)^{3/2}}e^{-\frac{x^2}{4(t-\tau)}}\left[1-\frac{x^2}{2(t-\tau)}\right]d\tau,$$

where from comes

$$-q_s = -\left(\frac{\partial T}{\partial x}\right)_0 = \frac{-1}{2\sqrt{\pi}}\int_0^t \frac{T_s(\tau)}{(t-\tau)^{3/2}}e^0 \cdot (1-0)d\tau = \frac{-1}{2\sqrt{\pi}}\int_0^t \frac{T_s(\tau)d\tau}{(t-\tau)^{3/2}}. \quad (32)$$

Comparing (32) with (27), one can see that they coincide because from (27) follows

$$-q_s = -\left(\frac{\partial T}{\partial x}\right)_0 = \frac{1}{\sqrt{\pi}}\frac{d}{dt}\int_0^t \frac{T_s(\tau)}{\sqrt{t-\tau}}d\tau = \frac{1}{\sqrt{\pi}}\left[\int_0^t \frac{-T_s(\tau)d\tau}{2(t-\tau)^{3/2}}+1\cdot\left(\frac{T_s(\tau)}{\sqrt{t-\tau}}\right)_{\tau=t}\right]=$$

$$= \frac{-1}{2\sqrt{\pi}}\int_0^t \frac{T_s(\tau)d\tau}{(t-\tau)^{3/2}}, \quad (33)$$

where $T_s(\tau) = 0$ until $t = \tau$. Thus, (32) and (33) coincide. The solutions obtained by the method of fractional differentials and the exact analytical solution of the boundary problem (24) completely coincide.

Now this heat flux at the boundary will be got ones more following to our algorithm. For this, first the temperature profile is computed through the Taylor series (17):

$$T = \left(\frac{\partial T}{\partial t}\right)_0 \Delta t + \frac{1}{2!}\left(\frac{\partial^2 T}{\partial t^2}\right)_0 (\Delta t)^2 + \frac{1}{3!}\left(\frac{\partial^3 T}{\partial t^3}\right)_0 (\Delta t)^3 + \ldots \quad (34)$$

where $T_0 = T(x,0) = 0$ according to (24). Replacing here $\frac{\partial T}{\partial t}$ through $\frac{\partial^2 T}{\partial x^2}$ in accordance with (24) results

$$T = \left(\frac{\partial^2 T}{\partial x^2}\right)_0 \Delta t + \frac{1}{2!}\left(\frac{\partial^4 T}{\partial x^4}\right)_0 (\Delta t)^2 + .... = \frac{\partial^2 T_s}{\partial x^2}\Delta t + \frac{1}{2!}\frac{\partial^4 T_s}{\partial x^4}(\Delta t)^2 + .... \quad (35)$$

Then the Taylor series expansion (34) with account of (35) can be rewritten up to an arbitrary order by $\Delta t$:

$$\Delta t: \quad T_1 = \sum_{n=1}^{\infty} \frac{1}{n!}\left(\frac{\partial^{2n} T}{\partial x^{2n}}\right)_0 (\Delta t)^n,$$

$$(2\Delta t): \quad T_2 = \sum_{n=1}^{\infty} \frac{1}{n!}\left(\frac{\partial^{2n} T}{\partial x^{2n}}\right)_0 (\Delta t)^n + \sum_{n=1}^{\infty} \frac{1}{n!}\left(\frac{\partial^{2n} T}{\partial x^{2n}}\right)_1 (\Delta t)^n, \quad (36)$$

$$m \cdot \Delta t: \quad T_m = T(t_m) = \sum_{n=1}^{\infty} \frac{1}{n!}\left[\left(\frac{\partial^{2n} T}{\partial x^{2n}}\right)_0 + \left(\frac{\partial^{2n} T}{\partial x^{2n}}\right)_1 + .... + \left(\frac{\partial^{2n} T}{\partial x^{2n}}\right)_m\right](\Delta t)^n$$

for the constant time step $\Delta t$ starting from $t = 0$ (actually, time step may be chosen variable, which does not affect the proposed method).

The Taylor series (36) represents the algorithm for numerical solution of the boundary problem (24) using approximations by time up to the desired accuracy. All needed for this is spatial derivatives by $x$ in the domain $x$. According to (35):

$$-q_s = -\left(\frac{\partial T}{\partial x}\right)_0 = -\left(\frac{\partial^3 T}{\partial x^3}\right)_0 t - \frac{1}{2!}\left(\frac{\partial^5 T}{\partial x^5}\right)_0 t^2 + .... \quad (37)$$

Now compute derivatives by $x$ using the exact analytical solution (30) of the outgoing boundary problem (24) and substitute into (37). Differentiating by $x$ yields from (30):

$$\frac{\partial^2 T}{\partial x^2} = \frac{1}{2\sqrt{\pi}}\int_0^t \frac{T_s(\tau)}{(t-\tau)^{3/2}} e^{-\frac{x^2}{4(t-\tau)}}\left[\frac{x}{2(t-\tau)}\left(\frac{x^2}{2(t-\tau)} - 3\right)\right]d\tau,$$

$$\frac{\partial^3 T}{\partial x^3} = \frac{1}{2\sqrt{\pi}}\int_0^t \frac{T_s(\tau)}{(t-\tau)^{3/2}} e^{-\frac{x^2}{4(t-\tau)}}\left[\frac{x^2}{4(t-\tau)^2}\left(3 - \frac{x^2}{2(t-\tau)}\right) + \frac{3}{4(t-\tau)}\left(\frac{x^2}{(t-\tau)} - 2\right)\right]d\tau,$$

$$\frac{\partial^4 T}{\partial x^4} = \frac{1}{2\sqrt{\pi}}\int_0^t \frac{T_s(\tau)}{(t-\tau)^{3/2}} e^{-\frac{x^2}{4(t-\tau)}}\left[\frac{x}{4(t-\tau)^2}\left(\frac{x^4}{4(t-\tau)^2} - \frac{5x^2}{t-\tau} + 15\right)\right]d\tau, \quad (38)$$

$$\frac{\partial^5 T}{\partial x^5} = \frac{1}{2\sqrt{\pi}}\int_0^t \frac{T_s(\tau)}{(t-\tau)^{3/2}} e^{-\frac{x^2}{4(t-\tau)}}\left\{\frac{-x^2}{8(t-\tau)^3}\left[\frac{x^4}{4(t-\tau)^2} - \frac{15x^2}{2(t-\tau)} + 45\right]\right\}d\tau,$$

$$\frac{\partial^6 T}{\partial x^6} = \frac{1}{2\sqrt{\pi}}\int_0^t \frac{T_s(\tau)}{(t-\tau)^{3/2}} e^{-\frac{x^2}{4(t-\tau)}}\left\{\frac{x}{4(t-\tau)^3}\left[\frac{x^6}{16(t-\tau)^3} - \frac{21x^4}{8(t-\tau)^2} + \frac{105x^2}{4(t-\tau)} - 45\right]\right\}d\tau,$$

$$\frac{\partial^7 T}{\partial x^7} = \frac{1}{2\sqrt{\pi}} \int_0^t \frac{T_s(\tau)}{(t-\tau)^{3/2}} e^{-\frac{x^2}{4(t-\tau)}} \left[ \frac{-x^8}{128(t-\tau)^7} + \frac{28x^6}{64(t-\tau)^6} - \frac{210x^4}{32(t-\tau)^5} + \right.$$

$$\left. + \frac{405x^2}{16(t-\tau)^4} - \frac{45}{4(t-\tau)^3} \right] d\tau,$$

The approximate numerical solution is given by recurrent formulae (36), therefore the analytic expression for derivatives (38) in equation (37) are satisfactory only for an initial small time step:

$$\left(\frac{\partial^3 T}{\partial x^3}\right)_0 = -\frac{3}{4\sqrt{\pi}} \int_0^t \frac{T_s(\tau)}{(t-\tau)^{5/2}} d\tau, \quad \left(\frac{\partial^5 T}{\partial x^5}\right)_0 = -\frac{3}{8\sqrt{\pi}} \int_0^t \frac{T_s(\tau)}{(t-\tau)^{5/2}} \left(\frac{1}{2} - \frac{5}{t-\tau}\right) d\tau. \quad (39)$$

Now the heat fluxes in a first and second order approaches are, respectively:

$$-q_s = -\left(\frac{\partial T}{\partial x}\right)_0 = -\frac{3t}{4\sqrt{\pi}} \int_0^t \frac{T_s(\tau)}{(t-\tau)^{5/2}} d\tau,$$

$$-q_s = -\left(\frac{\partial T}{\partial x}\right)_0 = -\frac{3t}{4\sqrt{\pi}} \int_0^t \frac{T_s(\tau)}{(t-\tau)^{5/2}} d\tau - \frac{3t^2}{16\sqrt{\pi}} \int_0^t \frac{T_s(\tau)}{(t-\tau)^{5/2}} \left(\frac{1}{2} - \frac{5}{t-\tau}\right) d\tau \quad (40)$$

where $t$ is small. To compare (40) with the exact solution, let us take also the third order approximation by $t$

$$-q_s = -\frac{3t}{4\sqrt{\pi}} \int_0^t \frac{T_s(\tau)}{(t-\tau)^{5/2}} d\tau + \frac{3t^2}{16\sqrt{\pi}} \int_0^t \frac{T_s(\tau)}{(t-\tau)^{5/2}} \left(\frac{1}{2} - \frac{5}{t-\tau}\right) d\tau + \frac{t^3}{64\sqrt{\pi}} \int_0^t \frac{T_s(\tau)}{(t-\tau)^{7/2}} \left(1 - \frac{42}{t-\tau}\right) d\tau. \quad (41)$$

The approximate solution (34) with accuracy of a third order terms by time is as follows:

$$T(x,t) = (1-\operatorname{sgn} x) T_s(t) \frac{x}{2\sqrt{\pi}} \int_0^t \frac{T_s(\tau)}{(t-\tau)^{3/2}} e^{-\frac{x^2}{4(t-\tau)}} \frac{1}{2} \left\{ \frac{t}{(t-\tau)} \left[ \frac{x^2}{2(t-\tau)} - 3 \right] + \frac{t^2}{4(t-\tau)^2} \cdot \right.$$

$$\left. \cdot \left[ \frac{x^4}{4(t-\tau)^2} - \frac{5x^2}{t-\tau} + 15 \right] + \frac{t^3}{12(t-\tau)^3} \left[ \frac{x^6}{16(t-\tau)^3} - \frac{21x^4}{8(t-\tau)^2} + \frac{105x^2}{4(t-\tau)} - 45 \right] \right\} d\tau, \quad (42)$$

where sgn($x$) is a signum function, $\operatorname{sgn} 0 = 0$, $\operatorname{sgn}(x) = 1$ by $x > 0$.

To estimate the approximate solution (42) obtained by our strategy and compare it with the exact analytical solution (27), the new variable $\xi = \frac{x}{2\sqrt{t-\tau}}$ is introduced in (42), which go to

$$\Delta T = \frac{2}{\sqrt{\pi}} \int_{\frac{x}{2\sqrt{t}}}^{\infty} T\left(t + \frac{x^2}{4\xi^2}\right) e^{-\xi^2} \frac{t\xi^2}{x^2} \left[ -4\xi^2 + 6 - 2\frac{t\xi^2}{x^2}\left(4\xi^4 - 2\xi^2 + 15\right) + \right.$$

$$\left. + \frac{8t^2 \xi^4}{3x^4} \left(4\xi^6 - 42\xi^4 + 105\xi^2 - 45\right) \right] d\xi, \quad (43)$$

where $\Delta T = T(x,t) - \tilde{T}(x,t)$, $\tilde{T}(x,t)$ is the approximate solution (42). In a first order by time from (43)

$$\Delta T = \left[1-(1-\operatorname{sgn} x)\right]T_s(t)+\ldots = (\operatorname{sgn} x)T_s(t)+\ldots = -\frac{x}{4\sqrt{t}}e^{-\frac{x^2}{4t}}, \qquad (44)$$

with a small exponential inaccuracy. By $x=0$, $\Delta T = 0$. If the multiplayer $\frac{2}{\sqrt{\pi}}T_s^*$ $\left(T_s^* = \max T_s\right)$ is got out of integral in the expression (43), then it goes to

$$\int_{\frac{x}{2\sqrt{t}}}^{\infty} e^{-\xi^2}\left[\frac{t\xi^2}{x^2}(6-4\xi^2)-2\frac{t^2\xi^4}{x^4}(4\xi^4-2\xi^2+15)-\frac{8t^3\xi^6}{3x^6}(4\xi^6-42\xi^4+105\xi^2-45)\right]d\xi \qquad (45)$$

Then in a first order accuracy by time $t$:

$$\operatorname{sgn}(x) - \int_0^{\frac{x}{2\sqrt{t}}} e^{-\xi^2}d\xi - \frac{x}{4\sqrt{t}}e^{-\frac{x^2}{4t}} = -\int_0^{\frac{x}{2\sqrt{t}}} e^{-\xi^2}d\xi - \frac{x}{4\sqrt{t}}e^{-\frac{x^2}{4t}}, \qquad (46)$$

where from by $x=0$, $t=0$ results exactly zero: 0-0-0=0 and 1-1-0=0.

Afterwards accounting the formulae (45) and (44), in a first order approach by time the following deficiency of the approximate solution is got:

$$\Delta T = \left(\int_0^{\infty} e^{-\xi^2}d\xi - \int_0^{\frac{x}{2\sqrt{t}}} e^{-\xi^2}d\xi\right)\operatorname{sgn}(x) - \frac{x}{4\sqrt{t}}e^{-\frac{x^2}{4t}}. \qquad (47)$$

In a second order approach by $t$, accordingly:

$$\Delta T = \left(\int_0^{\infty} e^{-\xi^2}d\xi - \int_0^{\frac{x}{2\sqrt{t}}} e^{-\xi^2}d\xi\right)\operatorname{sgn}(x) - \frac{x}{4\sqrt{t}}e^{-\frac{x^2}{4t}} - \frac{x}{16\sqrt{t}}\left(5+\frac{x^2}{2t}\right)e^{-\frac{x^2}{4t}}, \qquad (48)$$

where from follows that inaccuracy of the numerical solution is $\sim (x/\sqrt{t})e^{-\frac{x^2}{t}}$, which exceeds by small $t$ an order of $t$ dramatically, $(t\sqrt{t}/x)e^{\frac{x^2}{t}}$ times, exponentially. For instance, by $t=10^{-2}$ it is

$$\frac{10^{-2}\sqrt{10^{-2}}}{x}e^{\frac{x^2}{10^{-2}}} = \frac{10^{-3}}{x}e^{10^2 x^2}, \qquad \text{or} \qquad \frac{e^{100x^2}}{10^3 x},$$

Which is huge value, except small $x$, where $e^{\frac{x^2}{4t}} \approx 1 - \frac{x^2}{4t}$ and then $\Delta T = -\frac{x}{4\sqrt{t}}\left(1-\frac{x^2}{4t}\right)$, small value due to small $x$.

By *x*=0, $\Delta T = 0$, then by *t*=0, $\Delta T = 0$, and we got complete coincide with the exact solution. The second approach (48) decreases accuracy by *t* adding the term $(x^3/t\sqrt{t})e^{-\frac{x^2}{t}}$ of order $x^2/t$ comparing to the first order solution. Obviously, by fixed *x* and small *t* this can decrease accuracy, which is good in a first order approach. In (45) the following integrals were computed:

$$2\int_{\frac{x}{2\sqrt{t}}}^{\infty} \xi^2 e^{-\xi^2} d\xi = -\int_{\frac{x}{2\sqrt{t}}}^{\infty} \xi e^{-\xi^2} d(-\xi^2) = \int_{\frac{x}{2\sqrt{t}}}^{\infty} e^{-\xi^2} d\xi - \xi e^{-\xi^2} \Big|_{\frac{x}{2\sqrt{t}}}^{\infty} = \frac{x}{2\sqrt{t}} e^{-\frac{x^2}{4t}} + \int_{\frac{x}{2\sqrt{t}}}^{\infty} e^{-\xi^2} d\xi,$$

$$2\int_{\frac{x}{2\sqrt{t}}}^{\infty} \xi^4 e^{-\xi^2} d\xi = \int_{\frac{x}{2\sqrt{t}}}^{\infty} e^{-\xi^2} \xi^3 d\xi = \frac{x^3}{8t\sqrt{t}} e^{-\frac{x^2}{4t}} + \frac{3x}{4\sqrt{t}} e^{-\frac{x^2}{4t}} + \frac{3}{2}\int_{\frac{x}{2\sqrt{t}}}^{\infty} e^{-\xi^2} d\xi,$$

$$2\int_{\frac{x}{2\sqrt{t}}}^{\infty} e^{-\xi^2} \xi^6 d\xi = \int_{\frac{x}{2\sqrt{t}}}^{\infty} \xi^5 de^{-\xi^2} = -\left(\frac{x^5}{2^5 t^{5/2}} + \frac{5x^3}{16 t^{3/2}} + \frac{15x}{8\sqrt{t}}\right)e^{-\frac{x^2}{4t}} - \frac{15}{4}\int_{\frac{x}{2\sqrt{t}}}^{\infty} e^{-\xi^2} d\xi,$$

$$2\int_{\frac{x}{2\sqrt{t}}}^{\infty} e^{-\xi^2} \xi^8 d\xi = \frac{x}{\sqrt{t}}\left(\frac{x^6}{2^7 t^3} + \frac{7x^4}{2^6 t^2} + \frac{35x^2}{32t} + \frac{105}{16}\right)e^{-\frac{x^2}{4t}} + \frac{105}{8}\int_{\frac{x}{2\sqrt{t}}}^{\infty} e^{-\xi^2} d\xi,$$

$$\int_{\frac{x}{2\sqrt{t}}}^{\infty} 2\frac{t^2 \xi^4}{x^4} e^{-\xi^2}\left(4\xi^4 - 20\xi^2 + 15\right) d\xi = \frac{x}{2^4 \sqrt{t}}\left(5 + \frac{x^2}{2t}\right)e^{-\frac{x^2}{4t}}.$$

Surprisingly all terms of a negative power by *x* are mutually omitted, the same as the integrals, then only the exponent and the arguments $x/\sqrt{t}$, $x^2/t$ are kept. By *x*=0 the numerical solution by our algorithm coincides completely with the exact analytical solution in a first order accuracy by time. By fixed *x* the accuracy is high (deficiency decreases exponentially by time). By small *x* inaccuracy may grow, therefore it is important to choose right steps by *x* and *t* nearby the boundary *x*=0.

## 6. NUMERICAL RESULTS TO SUPPORT THE EFFICACY OF THE METHOD

Now the numerical solution by the method proposed here is applied to the non-stationary two-dimensional heat transfer problem for the Stokes flow around the sphere:

$$\left\{\frac{\partial}{\partial \tau} - \frac{\partial^2}{\partial \rho^2} - \left[\frac{2}{\rho} - \frac{\text{Pe}}{2}\cos\theta\left(1 - \frac{3}{2\rho} + \frac{1}{2\rho^3}\right)\right]\frac{\partial}{\partial \rho} - \frac{1}{\rho^2}\frac{\partial^2}{\partial \theta^2} - \left[\frac{ctg\theta}{\rho^2} + \frac{\text{Pe}\sin\theta}{2\rho}\left(1 - \frac{3}{4\rho} - \frac{1}{4\rho^3}\right)\right]\frac{\partial}{\partial \theta}\right\}T = 0,$$

(49)

$$1 < \rho < \infty, \quad 0 \leq \theta \leq \pi, \quad 0 < \tau < \infty; \quad T\big|_{\rho=1} = T_s(\theta,\tau); \quad T\big|_{\rho=\infty} = 0; \quad T\big|_{\tau=0} = 0,$$

with a symmetry conditions by $\theta$, where are $\rho = \frac{r}{R}$, $\tau = \frac{at}{R^2}$, $\text{Pe} = \frac{2UR}{a}$.

Numerical simulation performed for a number of a different boundary temperature distribution $T_s(\theta,\tau)$ has shown that the second and the third order by time solutions nearly coincide so that for this problem no higher than a second order by time is needed. A computation time among the first-fifth order by time does not reveal remarkable difference but by the sixth order by time it starts to grow substantially.

A few selected simulation results to support the efficacy of the method are given in the Table below:

A picture around the sphere is symmetrical, therefore the three points after $\theta = 4.3982$ are omitted just to save a place in the Table. For the initial-boundary value problem (49) there were not found any remarkable difference between computation results in the second and in the third order by time (for higher orders as well) despite diverse boundary conditions proven in the numerical simulation.

## 7. CONCLUSIONS

The examples considered here have clearly demonstrated an effectiveness and simplicity of the strategy proposed for the numerical solution of the non-stationary non-isothermal Navier-Stokes equations, as well as any other first order by time partial differential equations. An order of the equations by spatial variables has no matter for this algorithm. The strategy is based on application of a Taylor series by time for the computation of the solution sought by its temporal derivatives. These temporal derivatives of the functions are expressed from the outgoing equations through their spatial derivatives. The highest order temporal derivatives used in the Taylor series for computation of the approximate solution are got differentiating by time the outgoing PDEs. Only the one Poisson equation for the pressure distribution has to be solved numerically in case of the full Navier-Stokes equation array.

The method is applicable both for incompressible, as well as for compressible Navier-Stokes Equations. It allows comparably easily increasing the fidelity by time up to a fifth order reducing solution of the outgoing non-stationary problem to a solution of a consequence of the stationary problems at each temporal step.

The strategy of numerical solution of the boundary problems for the PDEs has been considered for a few diverse examples. An efficiency and simplicity of the method is achieved through a replacement of a solution of the finite-difference (or finite-element, etc.) equations in the known numerical methods with a simpler procedure of the spatial differentiation and further computation of the Taylor series. The strategy is going to be proved on different CFD problems to study all the pros and cons for its further implementation (e.g. restriction about the time step of the method, optimization between the time steps and the order of a numerical solution by time, etc.).

Table. Computation of a temperature distribution $T(r,\theta,\tau)$ in a Stokes flow around a sphere for a time dependent temperature on the sphere (r=1) $T_s(\theta,\tau) = |\cos\theta|\exp(100\tau)$.

| 1st order | $\tau$ =0.03 | | | | | | | |
|---|---|---|---|---|---|---|---|---|
| | $\theta$ = | | | | | | | |
| r = | 0.0000 | 0.6283 | 1.2566 | 1.8850 | 2.5133 | 3.1416 | 3.7699 | 4.3982 |
| 1 | 20.08554 | 16.24954 | 6.206772 | 6.206772 | 16.24954 | 20.08554 | 16.24954 | 6.206772 |
| 2 | 0.081860 | 0.115846 | 0.061748 | 0.060168 | 0.105751 | 0.178708 | 0.105759 | 0.059840 |
| 5 | 0.012018 | 0.012017 | 0.012013 | 0.012007 | 0.012003 | 0.012000 | 0.012003 | 0.012007 |
| 8 | 0.007506 | 0.007506 | 0.007504 | 0.007501 | 0.007499 | 0.007498 | 0.007499 | 0.007501 |
| 2nd | $\theta$ | | | | | | | |
| 1 | 20.08554 | 16.24954 | 6.206772 | 6.206772 | 16.24954 | 20.08554 | 16.24954 | 6.206772 |
| 2 | 0.081728 | 0.115994 | 0.061897 | 0.060317 | 0.105899 | 0.179139 | 0.105907 | 0.059987 |
| 5 | 0.012078 | 0.012077 | 0.012073 | 0.012067 | 0.012063 | 0.012060 | 0.012063 | 0.012067 |
| 8 | 0.007544 | 0.007543 | 0.007541 | 0.007539 | 0.007537 | 0.007536 | 0.007537 | 0.007539 |
| 3rd | $\theta$ | | | | | | | |
| 1 | 20.08554 | 16.24954 | 6.206772 | 6.206772 | 16.24954 | 20.08554 | 16.24954 | 6.206772 |
| 2 | 0.081728 | 0.115994 | 0.061897 | 0.060318 | 0.105900 | 0.179140 | 0.105907 | 0.059988 |
| 5 | 0.012078 | 0.012077 | 0.012073 | 0.012068 | 0.012063 | 0.012060 | 0.012063 | 0.012068 |
| 8 | 0.007544 | 0.007543 | 0.007541 | 0.007539 | 0.007537 | 0.007536 | 0.007537 | 0.007539 |